\documentclass[11pt,a4paper]{article}
\usepackage[T1]{fontenc}
\usepackage{lmodern,amsmath,amssymb,amsthm,graphicx}
\usepackage[margin=28mm]{geometry}
\usepackage[hidelinks]{hyperref}
\numberwithin{equation}{section}
\theoremstyle{plain}
\newtheorem{theorem}{Theorem}[section]
\theoremstyle{definition}
\newtheorem{remark}{Remark}[section]
\theoremstyle{plain}
\newtheorem{lemma}{Lemma}[section]
\newcommand{\R}{\mathbb R}
\newcommand{\T}{\mathbb T}
\newcommand{\norm}[1]{\lVert #1\rVert}

\newcommand{\dd}{\,d}
\newcommand{\proofpart}[1]{\par\smallskip\noindent\textbf{#1}\enspace\ignorespaces}
\title{A momentum criterion for wave breaking in the Camassa–Holm equation with linear dispersion}
\author{Song Liu\thanks{Email: lius37@mail2.sysu.edu.cn.}, Zhaonan Luo\thanks{Email: luozhn7@mail.sysu.edu.cn.} and Zhaoyang Yin\thanks{Email: mcsyzy@mail.sysu.edu.cn.}\\[4pt]\small School of Science, Shenzhen Campus of Sun Yat-sen University,\\\small Shenzhen 518107, China}
\date{}
\begin{document}
\maketitle
\begin{abstract}
We characterize wave breaking for the Camassa--Holm equation with
linear dispersion in terms of the initial momentum.
Let $u_0\in H^s(\mathbb R)$, $s>5/2$, be the initial velocity and
$m_0=u_0-u_0''+\lambda/2$ the initial momentum.
For $\lambda\ne0$, we prove that the strong solution undergoes wave
breaking in finite time if and only if $\lambda m_0(x_0)<0$ for some
$x_0\in\mathbb R$.
We also obtain corresponding characterizations for the periodic
Camassa--Holm equation and for the modified Camassa--Holm equation
on the real line and on the circle.
\end{abstract}
\noindent\textbf{Keywords:} Camassa--Holm equation; wave breaking;
initial momentum.
\tableofcontents

\setlength{\parskip}{3pt plus 0.5pt}
\section{Introduction}
\subsection{The equations and previous results}
The Camassa--Holm equation with linear dispersion is
\begin{equation}
  \left\{\begin{aligned} u_t-u_{txx}+3uu_x+\lambda u_x&=2u_xu_{xx}+uu_{xxx}, && t>0,\ x\in\R,\\ u(0,x)&=u_0(x), && x\in\R.\end{aligned}\right.
  \label{eq:1.1}
\end{equation}
Here $\lambda\in\R$ is the linear-dispersion parameter. We define the momentum by
\begin{equation}
  m=u-u_{xx}+\frac\lambda2,\qquad m_0=u_0-\partial_x^2u_0+\frac\lambda2.
  \label{eq:1.2}
\end{equation}
With $p(x)=e^{-|x|}/2$, the nonlocal form of \eqref{eq:1.1} is
\begin{equation}
  u_t+uu_x=-P_x,\qquad P=p*\left(u^2+\frac12u_x^2+\lambda u\right).
  \label{eq:1.3}
\end{equation}
We also consider \eqref{eq:1.1} on $\T$, with arbitrary $\lambda\in\R$. In this case, $p$ denotes the periodic Green function of $1-\partial_x^2$, and convolution is taken over one period. On either domain, $T^*\in(0,\infty]$ denotes the maximal existence time of the strong solution
\begin{equation}
  u\in C([0,T^*);H^s)\cap C^1([0,T^*);H^{s-1}),\qquad s>3/2.
  \label{eq:1.4}
\end{equation}
The characteristic flow is denoted by $q$:
\begin{equation}
  q_t(t,\xi)=u(t,q(t,\xi)),\qquad q(0,\xi)=\xi.
  \label{eq:1.5}
\end{equation}

The modified Camassa--Holm equation is considered on both $\R$ and $\T$. For $\lambda\ne0$, it takes the form
\begin{equation}
  \left\{\begin{aligned}\gamma_t+(u\gamma)_x&=\frac12\gamma^2+\lambda u,\\ u_{xx}-u&=\gamma_x+\frac{\gamma^2}{2\lambda},\\ \gamma(0,x)&=\gamma_0(x).\end{aligned}\right.
  \label{eq:1.6}
\end{equation}
The second equation determines the velocity $u$ uniquely in the relevant Sobolev space. Its initial value is denoted by $u_0$. The corresponding momentum is
\begin{equation}
  m=\frac{\lambda}{2}-\gamma_x-\frac{\gamma^2}{2\lambda},\qquad m_0=\frac{\lambda}{2}-\partial_x\gamma_0-\frac{\gamma_0^2}{2\lambda}.
  \label{eq:1.7}
\end{equation}
For $\gamma_0\in H^s$, $s>1/2$, we use $T^*$ for the maximal existence time in the class
\begin{equation}
  \gamma\in C([0,T^*);H^s)\cap C^1([0,T^*);H^{s-1}),\qquad s>1/2.
  \label{eq:1.8}
\end{equation}
The associated velocity belongs to $H^{s+1}$, $s>1/2$.

The CH equation emerged from the Hamiltonian description of shallow-water motion in the work of Camassa and Holm \cite{r5} and Camassa, Holm, and Hyman \cite{r6}; its integrable structure is connected with the earlier work of Fuchssteiner and Fokas \cite{r18}. A rigorous assessment of its water-wave approximation was given by Constantin and Lannes \cite{r12}. For the equation without linear dispersion, the theory includes local well-posedness in $H^s(\R)$, $s>3/2$ \cite{r26}, Besov-space well-posedness \cite{r16,r17,r34}, and norm inflation in $H^{3/2}(\R)$ \cite{r20}. Global strong solutions and wave breaking were studied in \cite{r9,r8,r31}. The peaked solitary waves are orbitally stable in $H^1(\R)$ \cite{r14}.

For nonzero dispersion, Constantin \cite{r7} studied local solutions and the scattering problem. The periodic theory was developed by Constantin and Escher \cite{r10} and Misio\l ek \cite{r32}. The local $H^s$ theory, $s>3/2$, is also recorded in Brandolese and Cortez \cite[Theorem 2.1]{r2}. For positive dispersion, Constantin and Strauss \cite{r15} proved orbital stability of smooth solitary waves. The inverse scattering transform and long-time asymptotics were studied in \cite{r11} and \cite{r1}, respectively.

The Miura correspondence with MOCH was developed by G\'orka and Reyes \cite{r19}. Luo, Qiao and Yin established the real-line local theory in Besov spaces \cite{r27} and the periodic local theory in integer-order $H^k$, $k\ge1$ \cite{r29}.

The conserved strong-solution energies are the CH $H^1$ norm and the MOCH $L^2$ norm; see \cite[Theorem 2.1]{r2}, \cite[Lemma 4.1]{r27} and \cite[Lemma 2.7]{r29}. We use $\norm{u}_{H^1}^2=\int(u^2+u_x^2)\dd x$ on the respective domain. For CH, global weak solutions were obtained in \cite{r13,r33}, and conservative and dissipative continuations were constructed in \cite{r3,r4,r23}. The conservative theory of Holden and Raynaud \cite{r24} also includes arbitrary linear dispersion and describes the terminal velocity and concentrated energy at wave breaking.

For MOCH, Luo, Qiao and Yin \cite[Theorem 1.3 and Corollary 1.4]{r28} constructed global conservative weak solutions for $L^2(\R)$ initial data. He and Yin \cite[Theorem 2.3 and Corollary 2.4]{r22} proved uniqueness by combining the energy balance with characteristic coordinates. These conservative theories supply the terminal weak limits and energy representations used in Section~4; the required statements will be quoted there without proof. For nonnegative dispersion, Kim, Kwon and Yoon \cite[Theorem 1.1]{r25} obtained $C^{3/5}$ breaking profiles for $H^5(\R)$ data satisfying their profile and localization assumptions.

For positive dispersion, Constantin \cite[Section 3(a)]{r7} discussed local solutions for $u_0\in H^k(\R)$, $k\in\mathbb N$, $k\ge3$, and showed that nonnegative shifted momentum is sufficient for global existence:
\begin{equation}
  \lambda>0,\quad m_0(x)\ge0\quad(x\in\R)\quad\Longrightarrow\quad T^*=\infty.
  \label{eq:1.9}
\end{equation}
He also established the wave-breaking alternative and the existence of smooth compactly supported initial profiles that break in finite time.

For CH with linear dispersion, Brandolese and Cortez \cite[Corollary 2.4]{r2} proved that, for $u_0\in H^s(\R)$, $s>3/2$,
\begin{equation}
  \partial_xu_0(x_0)<-\left|u_0(x_0)+\frac\lambda2\right|\quad\text{for some }x_0\in\R\quad\Longrightarrow\quad T^*<\infty.
  \label{eq:1.10}
\end{equation}
Their bounded-data result also covers periodic solutions \cite[Theorem 4.1 and the subsequent discussion]{r2}. They further obtained restrictions on global solutions \cite[Theorem 4.1 and Corollaries 4.3--4.4]{r2}. In particular, a global smooth CH solution satisfies
\begin{equation}
  u_x(t,x)\ge-\left|u(t,x)+\frac\lambda2\right|,
  \label{eq:1.11}
\end{equation}
and, in the periodic case,
\begin{equation}
  T^*=\infty,\quad u_0(x_0)=-\frac\lambda2\quad\Longrightarrow\quad u(t,x)\equiv-\frac\lambda2.
  \label{eq:1.12}
\end{equation}
These restrictions on global velocities are used below to obtain the initial-momentum characterizations.

For the equation without linear dispersion, McKean \cite[pp.~867--868]{r30} considered smooth periodic data and smooth compact perturbations of a constant, assuming finitely many proper changes of sign of the initial momentum. This assumption permits intervals on which the momentum vanishes. His breakdown criterion identifies the obstruction to global existence as a positive momentum value to the left of a negative one: writing $v_0$ for the initial velocity, there are $a<b$ such that
\[(v_0-\partial_x^2v_0)(a)>0,\qquad(v_0-\partial_x^2v_0)(b)<0.\]
In the periodic case, the corresponding global-existence criterion is that the initial momentum does not change sign. See also \cite{r31} for the later treatment.

The change of variables
\[v(t,x)=u\left(t,x-\frac{\lambda t}{2}\right)+\frac\lambda2\]
transforms \eqref{eq:1.1} into the equation without linear dispersion, and $v_0-\partial_x^2v_0=m_0$. For $u_0\in H^s(\R)$, $s>5/2$, we have $m_0(x)\to\lambda/2$ as $x\to\pm\infty$. Thus, when $\lambda\ne0$, the absence of the above sign configuration is equivalent to $\lambda m_0\ge0$. The real-line characterization below is proved for general $H^s(\R)$ data, $s>5/2$, without restricting the number of momentum sign changes; the transformed initial velocity belongs to $\lambda/2+H^s(\R)$. We also treat the periodic problem by the same characteristic argument.

For MOCH on $\R$, Luo, Qiao and Yin \cite[Theorem 4.3]{r27} proved finite-time blow-up for $\gamma_0\in H^s(\R)$, $s>1/2$, if
\begin{equation}
  \gamma_0(x_0)<-2\sqrt d\qquad\text{for some }x_0\in\R,
  \label{eq:1.13}
\end{equation}
where, with the positive constant $C$ from their convolution estimates,
\begin{equation}
  d=C\left(\norm{\gamma_0}_{L^2(\R)}+\frac{\norm{\gamma_0}_{L^2(\R)}^2}{2|\lambda|}\right)\left(C\norm{\gamma_0}_{L^2(\R)}+\frac{C\norm{\gamma_0}_{L^2(\R)}^2}{2|\lambda|}+|\lambda|\right).
  \label{eq:1.14}
\end{equation}
This condition compares a negative pointwise value of $\gamma_0$ with its conserved $L^2$ norm.

For periodic MOCH, Luo, Qiao and Yin \cite[Theorem 3.2 and its proof]{r29} established a small-data global result for $\gamma_0\in H^2(\T)$. Their estimates yield the sufficient condition
\begin{equation}
  \norm{\gamma_0}_{H^2(\T)}\le\sqrt{\frac23}(\sqrt2-1)|\lambda|,
  \label{eq:1.15}
\end{equation}
where this quoted bound uses the unit-period derivative norm
\begin{equation}
  \norm{\gamma_0}_{H^2(\T)}^2=\int_{\T}\!\bigl(\gamma_0^2+(\partial_x\gamma_0)^2+(\partial_x^2\gamma_0)^2\bigr)\dd x.
  \label{eq:1.16}
\end{equation}
For $\gamma_0\in H^k(\T)$, $k\in\mathbb N$, $k\ge1$, their blow-up theorem \cite[Theorem 4.2]{r29} states that
\begin{equation}
  \gamma_0(x_0)<-2\sqrt d\quad\text{for some }x_0\in\T\quad\Longrightarrow\quad T^*<\infty,
  \label{eq:1.17}
\end{equation}
where in this periodic result
\begin{align}
  d&=\left(3\norm{\gamma_0}_{L^2(\T)}+\frac{\norm{\gamma_0}_{L^2(\T)}^2}{2|\lambda|}\right)^2\notag\\
    &\quad+|\lambda|\left(\norm{\gamma_0}_{L^2(\T)}+\frac{3\norm{\gamma_0}_{L^2(\T)}^2}{2|\lambda|}\right).
  \label{eq:1.18}
\end{align}
Thus this work supplies both a small-data global theorem and a sufficiently negative-data blow-up theorem.

The classical breaking rate for the minimum slope in the dispersionless case is recalled in \cite[Theorem 1.1(c)]{r21}: for $u_0\in H^s(\R)$, $s\ge3$, and $T^*<\infty$,
\begin{equation}
  \lim_{t\uparrow T^*}(T^*-t)\inf_{x\in\R}u_x(t,x)=-2.
  \label{eq:1.19}
\end{equation}
McKean \cite[Section 4, pp.~872--873]{r30} also obtained the inverse-time slope behavior along breaking characteristics and discussed the limiting flow and compression of zero-momentum intervals under the assumptions stated above.
\subsection{Main results}
Our main results characterize global existence through the continuous initial momentum. We treat CH on the real line and on the circle for $u_0\in H^s$, $s>5/2$, and obtain the corresponding MOCH statements for $\gamma_0\in H^s$, $s>3/2$. The nonzero background fixes the admissible momentum sign on $\R$, whereas either sign is allowed on $\T$. The periodic CH condition recovers the classical one-sign criterion in the stated Sobolev framework. At first breaking, we identify the zero-momentum labels and relate their quadratic compression to the terminal energy defect.

\begin{theorem}[The real line]
Let $\lambda\ne0$ and $u_0\in H^s(\R)$ with $s>5/2$. Let $T^*$ be the maximal existence time of the strong solution of \eqref{eq:1.1}. Then
\begin{equation}
  T^*=\infty\quad\Longleftrightarrow\quad \lambda m_0(x)\ge0\quad\text{for every }x\in\R.
  \label{eq:1.20}
\end{equation}
\end{theorem}
\begin{remark}
The velocity has conserved $H^1$ energy and remains uniformly bounded:
\begin{align}
\norm{u(t)}_{H^1(\R)}&=\norm{u_0}_{H^1(\R)},\notag\\
     \norm{u(t)}_{L^\infty(\R)}&\le\frac1{\sqrt2}\norm{u_0}_{H^1(\R)},\qquad 0\le t<T^*. \label{eq:1.21}
\end{align}
In the global case,
\begin{equation}
  |u_x(t,x)|\le\left|u(t,x)+\frac\lambda2\right|,\qquad t\ge0.
  \label{eq:1.22}
\end{equation}
If $T^*<\infty$, then
\begin{equation}
  \liminf_{t\uparrow T^*}\inf_xu_x(t,x)=-\infty.
  \label{eq:1.23}
\end{equation}
\end{remark}
\begin{theorem}[The periodic problem]
Let $\lambda\in\R$ and $u_0\in H^s(\T)$ with $s>5/2$. Let $T^*$ be the maximal existence time of the periodic strong solution. Then
\begin{equation}
  T^*=\infty\quad\Longleftrightarrow\quad m_0\ge0\text{ on }\T\quad\text{or}\quad m_0\le0\text{ on }\T.
  \label{eq:1.24}
\end{equation}
\end{theorem}
\begin{remark}
The global derivative bound and wave-breaking conclusion in Remark 1.1 remain valid on $\T$. The energy and amplitude bounds are
\begin{align}
\norm{u(t)}_{H^1(\T)}&=\norm{u_0}_{H^1(\T)},\notag\\
     \norm{u(t)}_{L^\infty(\T)}&\le C\norm{u_0}_{H^1(\T)},\qquad 0\le t<T^*. \label{eq:1.25}
\end{align}
\end{remark}
\begin{theorem}[The modified Camassa--Holm equation]
Let $\lambda\ne0$ and $\gamma_0\in H^s$, $s>3/2$, on $\R$ or $\T$, and let $m_0$ be given by \eqref{eq:1.7}. Let $T^*$ be the maximal existence time of the corresponding strong solution of \eqref{eq:1.6}.
On the real line,
\begin{equation}
T^*=\infty\quad\Longleftrightarrow\quad\lambda m_0(x)\ge0\quad\text{for every }x\in\R.
  \label{eq:1.26}
\end{equation}
On the circle,
\begin{equation}
T^*=\infty\quad\Longleftrightarrow\quad m_0\ge0\text{ on }\T\quad\text{or}\quad m_0\le0\text{ on }\T.
  \label{eq:1.27}
\end{equation}
\end{theorem}
\begin{remark}
The associated velocity satisfies the conclusions of Remarks 1.1 and 1.2 on the respective domains. For $\gamma$,
\begin{align}
T^*=\infty&\quad\Longrightarrow\quad \sup_{t\ge0}\norm{\gamma(t)}_{L^\infty}<\infty,\notag\\
     T^*<\infty&\quad\Longrightarrow\quad \liminf_{t\uparrow T^*}\inf_x\gamma(t,x)=-\infty. \label{eq:1.28}
\end{align}
\end{remark}
\begin{theorem}[Terminal flow and compression]
Let $\lambda\ne0$ and $u_0\in H^s(\R)$, $s>5/2$, and suppose that $T^*<\infty$. The characteristic flow has a terminal value satisfying
\begin{gather}
  q(T^*,\cdot)\in C^1(\R), \label{eq:1.29} \\
  q_\xi(T^*,z)=0\quad\Longrightarrow\quad m_0(z)=0. \label{eq:1.30}
\end{gather}
For $Z=\{z\in\R:q_\xi(T^*,z)=0\}$, the compressed image satisfies
\begin{equation}
  \lim_{t\uparrow T^*}\frac{4|q(t,Z)|}{(T^*-t)^2}=\norm{u_0}_{H^1}^2-\norm{u(T^*)}_{H^1}^2,
  \label{eq:1.31}
\end{equation}
where $|q(t,Z)|$ denotes its Lebesgue measure and $u(T^*)$ is the terminal value supplied by the conservative theory \cite[Theorem 2.9]{r24}.
\end{theorem}
\begin{remark}[CH and MOCH energy defects]
Let $\lambda\ne0$, $\gamma_0\in H^s(\R)$, $s>3/2$, and $T^*<\infty$. Let $u$ be the associated CH velocity. With the weak terminal values supplied by the conservative theories \cite{r24,r28,r22},
\begin{equation}
  \norm{\gamma_0}_2^2-\norm{\gamma(T^*)}_2^2=\norm{u_0}_{H^1}^2-\norm{u(T^*)}_{H^1}^2.
  \label{eq:1.32}
\end{equation}
\end{remark}
\subsection{Motivation and strategy}
Our aim is a direct proof of the momentum characterizations in the stated Sobolev classes. The main step is to relate the initial momentum to the pointwise restriction satisfied by every global solution.

We describe the necessity argument on $\R$ for $\lambda>0$. The single-point criterion \eqref{eq:1.10}, applied at every time, gives
\begin{equation}
  T^*=\infty\quad\Longrightarrow\quad u_x(t,x)\ge-\left|u(t,x)+\frac\lambda2\right|\qquad(t\ge0,\ x\in\R).
  \label{eq:1.33}
\end{equation}
At a fixed time, $u(t,x)+\lambda/2\to\lambda/2>0$ as $x\to-\infty$. On any interval to the right of $a$ where $u(t,\cdot)+\lambda/2>0$,
\begin{align}
\partial_x\!\left[e^x\left(u(t,x)+\frac\lambda2\right)\right]&\ge0,\notag\\
     u(t,x)+\frac\lambda2&\ge e^{a-x}\left(u(t,a)+\frac\lambda2\right)>0. \label{eq:1.34}
\end{align}
This lower bound excludes a first zero, so $u+\lambda/2>0$ everywhere.
The characteristic flow and the key integrated quantity are
\begin{align}
q_t(t,\xi)&=u(t,q(t,\xi)),\qquad q(0,\xi)=\xi,\notag\\
     F(t,\xi)&=\int_0^t\frac{\dd r}{q_\xi(r,\xi)},\qquad t>0. \label{eq:1.35}
\end{align}
The positivity of $u+\lambda/2$ determines the sign in the integrated momentum identity:
\begin{equation}
  m_0(\xi)+2\partial_\xi\frac1{F(t,\xi)}\ge0.
  \label{eq:1.36}
\end{equation}
Energy conservation bounds the velocity and the length of each transported interval:
\begin{align}
\norm{u(t)}_{L^\infty(\R)}&\le\frac1{\sqrt2}\norm{u_0}_{H^1(\R)},\notag\\
     0<q(t,b)-q(t,a)&\le b-a+\sqrt2\norm{u_0}_{H^1(\R)}t. \label{eq:1.37}
\end{align}
Consequently, an integral Cauchy--Schwarz estimate gives, for every fixed label interval $[a,b]$ with $a<b$,
\begin{equation}
  0<\int_a^b\frac{\dd\xi}{F(t,\xi)}\le\left(\int_0^t\frac{\dd r}{b-a+\sqrt2\norm{u_0}_{H^1(\R)}r}\right)^{-1}\longrightarrow0\qquad(t\to\infty).
  \label{eq:1.38}
\end{equation}

If $m_0$ is negative somewhere, continuity supplies a fixed interval $[a,b]$ on which $m_0<0$. Two integrations of the momentum inequality yield
\begin{equation}
  \int_a^b\frac{\dd\xi}{F(t,\xi)}\ge-\frac12\int_a^b(b-y)m_0(y)\dd y>0,\qquad t>0.
  \label{eq:1.39}
\end{equation}
The last lower bound is independent of time and contradicts the preceding limit. The sufficient condition follows from preservation of the momentum sign and the resulting bound on $u_x$.

Reflection handles negative dispersion. In the periodic problem, the pointwise restriction first excludes sign changes of $u+\lambda/2$; its conserved mean then selects the sign. Preservation of the period gives the required integral decay. For MOCH, the Miura correspondence and the conserved $L^2$ norm bound $\gamma-u_x$ uniformly. This proves equality of the maximal existence times and transfers the CH characterizations to MOCH.

At a finite maximal time, the quantity $\rho=\sqrt{q_\xi}$ satisfies a linear second-order equation in time with bounded coefficient. Its label derivative remains bounded for $u_0\in H^s$, $s>5/2$, so the terminal momentum identity gives
\begin{equation}
  \rho(T^*,z)=\rho_\xi(T^*,z)=0\quad\Longrightarrow\quad m_0(z)=0.
  \label{eq:1.40}
\end{equation}
In the initial label, the energy is
\begin{equation}
  \int_\R\!\bigl(u(t,q(t,\xi))^2\rho(t,\xi)^2+4\rho_t(t,\xi)^2\bigr)\dd\xi=\norm{u_0}_{H^1}^2.
  \label{eq:1.41}
\end{equation}
The conservative energy representation recalled above identifies the energy carried by labels with $q_\xi(T^*,\xi)=0$. Comparing the terminal weighted products $(\gamma\circ q)\rho$ and $2\rho_t$ on these labels then gives equality of the CH and MOCH energy defects.

\noindent The remaining sections are arranged as follows.

\noindent\textbf{Section 2.} We recall the local theory, conservation laws, and pointwise criterion. We derive the finite-time characterization and the flow and momentum identities, and establish the common maximal lifespan of CH and MOCH.

\noindent\textbf{Section 3.} We prove the momentum characterizations on the real line and on the circle, transfer them to MOCH, and establish the stated properties of global solutions and wave breaking.

\noindent\textbf{Section 4.} We establish the terminal flow, identify the zero-momentum breaking labels, and relate their quadratic compression to the terminal energy defect. We then compare the CH and MOCH energy defects, using the conservative results recalled in Section~1.1.

\noindent\textbf{Section 5.} We explain two numerical diagnostics for compression and mesh resolution.

\section{Preliminaries}
The CH statements hold on $\R$ and $\T$ for $\lambda\in\R$. We omit the domain when it is clear and use periodic lifts in integrals over $\R$.
\subsection{Local solutions}
For CH on $\R$, the local theory and one-sided continuation criterion are given by Brandolese and Cortez \cite[Theorem 2.1(1) and (2)(i)]{r2}. Persistence of higher regularity is discussed immediately after that theorem.

For the periodic problem, see Constantin and Escher \cite{r10} and Misio\l ek \cite{r32}. The $H^s$ local theory, $s>3/2$, is also recorded in \cite[Section 4, paragraph preceding Corollary 4.4]{r2}.
\begin{lemma}[Local solutions for CH]
Let $\lambda\in\R$ and $u_0\in H^s$, $s>3/2$. There exists a unique maximal strong solution of \eqref{eq:1.3} such that
\begin{equation}
  u\in C([0,T^*);H^s)\cap C^1([0,T^*);H^{s-1}),\qquad s>3/2.
  \label{eq:2.1}
\end{equation}
The solution depends continuously on the initial data on every compact subinterval of $[0,T^*)$, and higher Sobolev regularity persists on the same maximal interval. If $T^*<\infty$, then
\begin{equation}
  \liminf_{t\uparrow T^*}\inf_xu_x(t,x)=-\infty.
  \label{eq:2.2}
\end{equation}
\end{lemma}
For MOCH on $\R$, we use Luo, Qiao and Yin \cite[Theorem 3.2 and Lemma 4.2]{r27}, with
\begin{equation}
  H^s=B^s_{2,2},\qquad s>1/2.
  \label{eq:2.3}
\end{equation}
On $\T$, the same local theory and continuation criterion for $s>1/2$ follow by carrying out the arguments of \cite[Theorem 3.2 and Lemma 4.2]{r27} with periodic Littlewood--Paley blocks and Fourier multipliers; see \cite[Section 2]{r29} for the periodic estimates. The integer-order periodic results are \cite[Lemmas 2.6 and 4.1]{r29}.
\begin{lemma}[Local solutions for MOCH]
Let $\lambda\ne0$ and $\gamma_0\in H^s$, $s>1/2$, on $\R$ or $\T$. There exists a unique maximal solution of \eqref{eq:1.6} such that
\begin{equation}
  \gamma\in C([0,T^*);H^s)\cap C^1([0,T^*);H^{s-1}),\qquad s>1/2.
  \label{eq:2.4}
\end{equation}
The solution depends continuously on the initial data on compact subintervals of $[0,T^*)$, and higher regularity persists. If $T^*<\infty$, then
\begin{equation}
  \limsup_{t\uparrow T^*}\norm{\gamma(t)}_{L^\infty}=\infty.
  \label{eq:2.5}
\end{equation}
\end{lemma}
\subsection{Energy conservation}
For CH energy conservation, see Brandolese and Cortez \cite[Theorem 2.1(1)]{r2} on $\R$. For periodic energy and mean conservation, see \cite[Section 4, paragraph preceding Corollary 4.4]{r2}.
\begin{lemma}[CH energy]
Let $u_0\in H^s$, $s>3/2$, and let $u$ be the solution in Lemma 2.1. For $0\le t<T^*$,
\begin{equation}
  \int\bigl(u(t,x)^2+u_x(t,x)^2\bigr)\dd x=\int\bigl(u_0(x)^2+(\partial_xu_0(x))^2\bigr)\dd x.
  \label{eq:2.6}
\end{equation}
The Sobolev bounds are
\begin{equation}
  \norm{u(t)}_{L^\infty(\R)}\le\frac1{\sqrt2}\norm{u_0}_{H^1(\R)},\qquad \norm{u(t)}_{L^\infty(\T)}\le C\norm{u_0}_{H^1(\T)},
  \label{eq:2.7}
\end{equation}
on the respective domains. In the periodic case, the mean is also conserved:
\begin{equation}
  \int_\T u(t,x)\dd x=\int_\T u_0(x)\dd x.
  \label{eq:2.8}
\end{equation}
\end{lemma}
For MOCH, see Luo, Qiao and Yin \cite[Lemma 4.1]{r27} on $\R$ and \cite[Lemma 2.7]{r29} on $\T$. The fractional periodic statement follows by smooth approximation and continuous dependence.
\begin{lemma}[MOCH energy]
Let $\lambda\ne0$ and $\gamma_0\in H^s$, $s>1/2$, on $\R$ or $\T$. The solution in Lemma 2.2 satisfies
\begin{equation}
  \norm{\gamma(t)}_{L^2}=\norm{\gamma_0}_{L^2},\qquad 0\le t<T^*.
  \label{eq:2.9}
\end{equation}
\end{lemma}
\subsection{The pointwise criterion and finite-time characterization}
The criterion is Brandolese and Cortez \cite[Corollary 2.4]{r2}; the lifespan bound is \cite[Theorem 2.2, formula (2.5)]{r2}.

Both statements extend to periodic solutions by \cite[Theorem 4.1 and the paragraph following it]{r2}. Time translation allows any initial time $t_0$.
\begin{lemma}[A single-point blow-up criterion]
Let $u_0\in H^s$, $s>3/2$, and let $u$ be the maximal CH solution. If there exist $t_0\in[0,T^*)$ and $x_0$ such that
\begin{equation}
  u_x(t_0,x_0)<-\left|u(t_0,x_0)+\frac\lambda2\right|,
  \label{eq:2.10}
\end{equation}
then $T^*<\infty$ and
\begin{equation}
  T^*-t_0\le\frac2{\sqrt{u_x(t_0,x_0)^2-(u(t_0,x_0)+\lambda/2)^2}}.
  \label{eq:2.11}
\end{equation}
\end{lemma}
\smallskip\noindent\textbf{Finite-time characterization.}

For the CH solution with $u_0\in H^s$, $s>3/2$, one has
\begin{equation}
  T^*<\infty\quad\Longleftrightarrow\quad\exists\,t_0\in[0,T^*),\ x_0:\quad u_x(t_0,x_0)<-\left|u(t_0,x_0)+\frac\lambda2\right|.
  \label{eq:2.12}
\end{equation}
\begin{proof}
The reverse implication is Lemma 2.5. If $T^*<\infty$, Lemma 2.3 and \eqref{eq:2.2} give
\begin{equation}
  \sup_{0\le t<T^*}\norm{u(t)}_\infty<\infty,\qquad\liminf_{t\uparrow T^*}\inf_xu_x(t,x)=-\infty.
  \label{eq:2.13}
\end{equation}
Hence there exist $t_0<T^*$ and $x_0$ such that
\begin{align}
u_x(t_0,x_0)&<-\sup_{0\le t<T^*}\norm{u(t)}_{L^\infty}-\frac{|\lambda|}{2}\notag\\
     &\le-\left|u(t_0,x_0)+\frac\lambda2\right|. \label{eq:2.14}
\end{align}
This proves the forward implication without assuming attainment of the spatial infimum.
\end{proof}
\subsection{The flow and momentum identities}
\begin{lemma}[The characteristic flow]
Let $u_0\in H^s$, $s>3/2$. The flow
\begin{equation}
  q_t(t,\xi)=u(t,q(t,\xi)),\qquad q(0,\xi)=\xi
  \label{eq:2.15}
\end{equation}
is an increasing $C^1$ diffeomorphism in the label for each $t<T^*$; in the periodic case we use its lift. Moreover,
\begin{equation}
q_{\xi t}(t,\xi)=u_x(t,q(t,\xi))q_\xi(t,\xi),\qquad q_\xi(t,\xi)=\exp\!\left(\int_0^tu_x(r,q(r,\xi))\dd r\right)>0.
  \label{eq:2.16}
\end{equation}
\end{lemma}
\begin{proof}
Fix $\tau<T^*$. Since $u_0\in H^s$, $s>3/2$, Sobolev embedding gives
\begin{equation}
  \sup_{0\le t\le\tau}\bigl(\norm{u(t)}_\infty+\norm{u_x(t)}_\infty\bigr)<\infty.
  \label{eq:2.17}
\end{equation}
The characteristic equation therefore has a unique solution on $[0,\tau]$. Differentiation of its integral form yields
\begin{align}
q(t,\xi)&=\xi+\int_0^t u(r,q(r,\xi))\dd r,\notag\\
     q_\xi(t,\xi)&=1+\int_0^t u_x(r,q(r,\xi))q_\xi(r,\xi)\dd r\notag\\
     &=\exp\!\left(\int_0^t u_x(r,q(r,\xi))\dd r\right)\ge\exp\!\left(-\int_0^t\norm{u_x(r)}_\infty\dd r\right)>0. \label{eq:2.18}
\end{align}
Moreover,
\begin{align}
|q(t,\xi)-\xi|&\le\int_0^t\norm{u(r)}_\infty\dd r,\notag\\
 \lim_{\xi\to-\infty}q(t,\xi)&=-\infty,\qquad\lim_{\xi\to+\infty}q(t,\xi)=+\infty. \label{eq:2.19}
\end{align}
Thus $q(t,\cdot)$ is increasing and onto, and its inverse is $C^1$. The same argument applies to the periodic lift.
\end{proof}
\begin{lemma}[Momentum transport and reconstruction]
Let $u_0\in H^s$, $s>5/2$, and set $m=u-u_{xx}+\lambda/2$. Then
\begin{equation}
  m(t,q(t,\xi))q_\xi(t,\xi)^2=m_0(\xi).
  \label{eq:2.20}
\end{equation}
Moreover,
\begin{align}
u(t,x)+\lambda/2-u_x(t,x)&=\int_{-\infty}^x e^{y-x}m(t,y)\dd y,\notag\\
     u(t,x)+\lambda/2+u_x(t,x)&=\int_x^\infty e^{x-y}m(t,y)\dd y. \label{eq:2.21}
\end{align}
The integrals use the periodic lift on $\T$ and require no $L^1$ assumption on $m$.
\end{lemma}
\begin{proof}
For smooth solutions, \eqref{eq:1.1} and \eqref{eq:2.16} give
\begin{align}
m_t+um_x+2u_xm&=u_t-u_{txx}+3uu_x+\lambda u_x-uu_{xxx}-2u_xu_{xx}=0,\notag\\
 \frac{d}{dt}\bigl(m(t,q(t,\xi))q_\xi(t,\xi)^2\bigr)&=(m_t+um_x+2u_xm)(t,q(t,\xi))q_\xi(t,\xi)^2=0. \label{eq:2.22}
\end{align}

For $u_0\in H^s$, $s>5/2$, take smooth $u_0^{(n)}\to u_0$ in $H^s$. Lemma 2.1 and $H^s\hookrightarrow C_b^2$ imply, for every $\tau<T^*$,
\begin{align}
u^{(n)}&\longrightarrow u &&\text{in } C([0,\tau];H^s),\qquad s>5/2,\notag\\
 u^{(n)},u_x^{(n)},u_{xx}^{(n)}&\longrightarrow u,u_x,u_{xx} &&\text{uniformly on } [0,\tau]\times\R. \label{eq:2.23}
\end{align}
The flow integral equation and Gronwall's inequality give
\begin{align}
\sup_\xi|q^{(n)}(t,\xi)-q(t,\xi)|&\le\int_0^t\norm{u^{(n)}(r)-u(r)}_\infty\dd r\notag\\
 &\quad+\int_0^t\norm{u_x(r)}_\infty\sup_\xi|q^{(n)}(r,\xi)-q(r,\xi)|\dd r, \label{eq:2.24}
\end{align}
Thus
\begin{equation}
\sup_{0\le t\le\tau,\,\xi\in\R}|q^{(n)}(t,\xi)-q(t,\xi)|\longrightarrow0.
  \label{eq:2.25}
\end{equation}
By \eqref{eq:2.16},
\begin{align}
\left|\log q_\xi^{(n)}(t,\xi)-\log q_\xi(t,\xi)\right|&\le\int_0^t\norm{u_x^{(n)}(r)-u_x(r)}_\infty\dd r\notag\\
 &\quad+\int_0^t\norm{u_{xx}(r)}_\infty\sup_\xi|q^{(n)}(r,\xi)-q(r,\xi)|\dd r\longrightarrow0. \label{eq:2.26}
\end{align}
The convergence is uniform on $[0,\tau]\times\R$. Since $m-\lambda/2\in C([0,\tau];H^{s-2})$, $s>5/2$, the functions $m(t,\cdot)$ are uniformly equicontinuous for $t\in[0,\tau]$. Thus
\begin{align}
m_0(\xi)&=\lim_{n\to\infty}m_0^{(n)}(\xi)\notag\\
     &=\lim_{n\to\infty}m^{(n)}(t,q^{(n)}(t,\xi))\bigl(q_\xi^{(n)}(t,\xi)\bigr)^2=m(t,q(t,\xi))q_\xi(t,\xi)^2. \label{eq:2.27}
\end{align}
For the periodic problem these calculations use the periodic lifts.

The reconstruction formulas follow from the exact derivatives
\begin{align}
\partial_y\!\left[e^y\left(u(t,y)-u_x(t,y)+\frac\lambda2\right)\right]&=e^ym(t,y),\notag\\
 \partial_y\!\left[e^{-y}\left(u(t,y)+u_x(t,y)+\frac\lambda2\right)\right]&=-e^{-y}m(t,y). \label{eq:2.28}
\end{align}
Since $u$ and $u_x$ are bounded,
\begin{align}
\lim_{y\to-\infty}e^y\left(u(t,y)-u_x(t,y)+\frac\lambda2\right)&=0,\notag\\
 \lim_{y\to+\infty}e^{-y}\left(u(t,y)+u_x(t,y)+\frac\lambda2\right)&=0. \label{eq:2.29}
\end{align}
Integrating \eqref{eq:2.28}, we obtain
\begin{align}
\int_{-\infty}^x e^{y-x}m(t,y)\dd y&=e^{-x}\left[e^y\left(u(t,y)-u_x(t,y)+\frac\lambda2\right)\right]_{-\infty}^x=u(t,x)-u_x(t,x)+\frac\lambda2,\notag\\
 \int_x^\infty e^{x-y}m(t,y)\dd y&=-e^x\left[e^{-y}\left(u(t,y)+u_x(t,y)+\frac\lambda2\right)\right]_x^\infty=u(t,x)+u_x(t,x)+\frac\lambda2. \label{eq:2.30}
\end{align}
\end{proof}
\subsection{An integrated flow identity}
Define
\begin{equation}
  F(t,\xi)=\int_0^t\frac{\dd r}{q_\xi(r,\xi)},\qquad 0<t<T^*.
  \label{eq:2.31}
\end{equation}
\begin{lemma}[Integrated flow estimates]
If $u_0\in H^s$, $s>3/2$, then for every $a<b$,
\begin{equation}
  \int_a^b\frac{\dd\xi}{F(t,\xi)}\le\left(\int_0^t\frac{\dd r}{q(r,b)-q(r,a)}\right)^{-1}.
  \label{eq:2.32}
\end{equation}
If, in addition, $u_0\in H^s$ with $s>5/2$, then $F(t,\cdot)\in C^1$ and
\begin{equation}
m_0(\xi)+2\partial_\xi\frac1{F(t,\xi)}=\frac2{F(t,\xi)^2}\int_0^t\!\left(u(r,q(r,\xi))+\frac\lambda2\right)q_\xi(r,\xi)\bigl(F(t,\xi)-F(r,\xi)\bigr)\dd r.
  \label{eq:2.33}
\end{equation}
\end{lemma}
\begin{proof}
For each fixed label, Cauchy--Schwarz gives
\begin{equation}
  \left(\int_0^t\frac{\dd r}{q(r,b)-q(r,a)}\right)^2\le F(t,\xi)\int_0^t\frac{q_\xi(r,\xi)}{(q(r,b)-q(r,a))^2}\dd r.
  \label{eq:2.34}
\end{equation}
Dividing by $F(t,\xi)$ and integrating over $[a,b]$, we obtain
\begin{align}
&\left(\int_0^t\frac{\dd r}{q(r,b)-q(r,a)}\right)^2\int_a^b\frac{\dd\xi}{F(t,\xi)}\notag\\
     &\qquad\le\int_0^t\frac{\int_a^b q_\xi(r,\xi)\dd\xi}{(q(r,b)-q(r,a))^2}\dd r=\int_0^t\frac{\dd r}{q(r,b)-q(r,a)}. \label{eq:2.35}
\end{align}
All factors are positive, so division proves \eqref{eq:2.32}.

For $u_0\in H^s$, $s>5/2$,
\begin{equation}
  u\in C([0,T^*);H^s)\subset C([0,T^*);C_b^2),\qquad q(t,\cdot)\in C^2,\qquad s>5/2.
  \label{eq:2.36}
\end{equation}
By \eqref{eq:2.16} and \eqref{eq:2.20},
\begin{align}
\partial_t\partial_\xi\log q_\xi(t,\xi)&=u_{xx}(t,q(t,\xi))q_\xi(t,\xi)\notag\\
     &=\left(u(t,q(t,\xi))+\frac\lambda2\right)q_\xi(t,\xi)-\frac{m_0(\xi)}{q_\xi(t,\xi)}. \label{eq:2.37}
\end{align}
Since $q_\xi(0,\xi)=1$, integration gives
\begin{equation}
  \partial_\xi\log q_\xi(t,\xi)=\int_0^t\left(u(r,q(r,\xi))+\frac\lambda2\right)q_\xi(r,\xi)\dd r-m_0(\xi)F(t,\xi).
  \label{eq:2.38}
\end{equation}

Consequently, differentiating \eqref{eq:2.31} gives
\begin{align}
F_\xi(t,\xi)&=-\int_0^t\frac{\partial_\xi\log q_\xi(r,\xi)}{q_\xi(r,\xi)}\dd r\notag\\
 &=m_0(\xi)\int_0^t\frac{F(r,\xi)}{q_\xi(r,\xi)}\dd r\notag\\
 &\quad-\int_0^t\frac{1}{q_\xi(r,\xi)}\int_0^r\left(u(\sigma,q(\sigma,\xi))+\frac\lambda2\right)q_\xi(\sigma,\xi)\dd\sigma\dd r\notag\\
 &=\frac12m_0(\xi)F(t,\xi)^2\notag\\
 &\quad-\int_0^t\left(u(\sigma,q(\sigma,\xi))+\frac\lambda2\right)q_\xi(\sigma,\xi)\int_\sigma^t\frac{\dd r}{q_\xi(r,\xi)}\dd\sigma\notag\\
 &=\frac12m_0(\xi)F(t,\xi)^2\notag\\
 &\quad-\int_0^t\left(u(\sigma,q(\sigma,\xi))+\frac\lambda2\right)q_\xi(\sigma,\xi)\bigl(F(t,\xi)-F(\sigma,\xi)\bigr)\dd\sigma. \label{eq:2.39}
\end{align}
On each compact time and label interval, the integrands are continuous and $\inf q_\xi>0$. Thus differentiation and Fubini's theorem apply. Moreover,
\begin{equation}
  \partial_rF=\frac1{q_\xi},\qquad F(0,\xi)=0,\qquad \int_0^t\frac{F(r,\xi)}{q_\xi(r,\xi)}\dd r=\frac12F(t,\xi)^2.
  \label{eq:2.40}
\end{equation}
Multiplication of \eqref{eq:2.39} by $-2/F(t,\xi)^2$ proves \eqref{eq:2.33}.
\end{proof}
\subsection{The Miura correspondence and maximal lifespan}
The smooth MOCH-to-CH correspondence is established by G\'orka and Reyes \cite[Section 5.2, formulas (5.7)--(5.12)]{r19}. We pass this correspondence to the Sobolev class and show that the maximal existence times agree.
\begin{lemma}[The Miura correspondence and common lifespan]
Let $\lambda\ne0$ and $\gamma_0\in H^s$, $s>3/2$, on $\R$ or $\T$. For the strong MOCH solution, with $p$ the Green function on that domain,
\begin{align}
u&=-p_x*\gamma-\frac1{2\lambda}p*\gamma^2,\notag\\
     u_x&=\gamma-p*\gamma-\frac1{2\lambda}p_x*\gamma^2,\notag\\
     m&=\frac\lambda2-\gamma_x-\frac{\gamma^2}{2\lambda}. \label{eq:2.41}
\end{align}
The velocity belongs to
\begin{equation}
  u\in C([0,T^*);H^{s+1})\cap C^1([0,T^*);H^s),\qquad s>3/2,
  \label{eq:2.42}
\end{equation}
and is the CH solution with initial value
\begin{equation}
  u_0=-p_x*\gamma_0-\frac1{2\lambda}p*\gamma_0^2.
  \label{eq:2.43}
\end{equation}
The two maximal existence times coincide. On the common interval,
\begin{equation}
  \norm{\gamma(t)-u_x(t)}_{L^\infty}\le\norm{p}_{L^2}\norm{\gamma_0}_{L^2}+\frac{\norm{p_x}_{L^\infty}}{2|\lambda|}\norm{\gamma_0}_{L^2}^2,\qquad 0\le t<T^*.
  \label{eq:2.44}
\end{equation}
On $\R$, $\norm{p}_2=1/2$ and $\norm{p_x}_\infty=1/2$.
\end{lemma}
\begin{proof}
Let $T_\gamma^*$ denote the MOCH lifespan. Fix $\tau<T_\gamma^*$ and take smooth $\gamma_0^{(n)}\to\gamma_0$ in $H^s$, $s>3/2$. By Lemma 2.2,
\begin{equation}
  \gamma^{(n)}\longrightarrow\gamma\quad\text{in }C([0,\tau];H^s),\qquad s>3/2.
  \label{eq:2.45}
\end{equation}
The operators $p_x*$ and $p*$ gain one and two Sobolev derivatives, respectively. The algebra property of $H^s$, $s>3/2$, therefore yields
\begin{align}
\norm{u^{(n)}-u}_{H^{s+1}}&\le C\norm{\gamma^{(n)}-\gamma}_{H^s}+\frac C{|\lambda|}\norm{(\gamma^{(n)})^2-\gamma^2}_{H^s}\notag\\
     &\le C\left(1+\frac{\norm{\gamma^{(n)}}_{H^s}+\norm{\gamma}_{H^s}}{|\lambda|}\right)\norm{\gamma^{(n)}-\gamma}_{H^s},\qquad s>3/2. \label{eq:2.46}
\end{align}
Hence
\begin{equation}
  u^{(n)}\longrightarrow u\quad\text{in }C([0,\tau];H^{s+1}),\qquad s>3/2.
  \label{eq:2.47}
\end{equation}
The three identities in the lemma pass to the limit. By the cited smooth Miura correspondence, each $u^{(n)}$ satisfies
\begin{equation}
  u^{(n)}(t)=u_0^{(n)}-\int_0^t\!\left[u^{(n)}u_x^{(n)}+p_x*\left((u^{(n)})^2+\frac12(u_x^{(n)})^2+\lambda u^{(n)}\right)\right](r)\dd r.
  \label{eq:2.48}
\end{equation}
The integrands converge in $C([0,\tau];H^s)$, $s>3/2$. Hence
\begin{equation}
  u(t)=u_0-\int_0^t\!\left[uu_x+p_x*(u^2+\frac12u_x^2+\lambda u)\right](r)\dd r.
  \label{eq:2.49}
\end{equation}
Thus \eqref{eq:2.42} holds. By Lemma 2.1, $u$ is the CH solution with initial value \eqref{eq:2.43}. Writing its lifespan as $T_u^*$, we have
\begin{equation}
  T_\gamma^*\le T_u^*.
  \label{eq:2.50}
\end{equation}

To compare lifespans, apply Young's inequality to the second identity in \eqref{eq:2.41} and use the conservation law \eqref{eq:2.9}:
\begin{align}
\norm{\gamma(t)-u_x(t)}_\infty&\le\norm p_2\norm{\gamma(t)}_2+\frac{\norm{p_x}_\infty}{2|\lambda|}\norm{\gamma(t)}_2^2\notag\\
     &=\norm p_2\norm{\gamma_0}_2+\frac{\norm{p_x}_\infty}{2|\lambda|}\norm{\gamma_0}_2^2. \label{eq:2.51}
\end{align}
If $T_\gamma^*<T_u^*$, the CH solution exists in $H^{s+1}$, $s>3/2$, on the compact interval $[0,T_\gamma^*]$. It follows that
\begin{equation}
  \sup_{0\le t<T_\gamma^*}\norm{\gamma(t)}_\infty\le\sup_{0\le t\le T_\gamma^*}\norm{u_x(t)}_\infty+\norm p_2\norm{\gamma_0}_2+\frac{\norm{p_x}_\infty}{2|\lambda|}\norm{\gamma_0}_2^2<\infty,
  \label{eq:2.52}
\end{equation}
This contradicts the MOCH continuation criterion \eqref{eq:2.5}, so $T_\gamma^*=T_u^*$.

At the initial time, the last identity in \eqref{eq:2.41} gives
\begin{equation}
  \lambda m_0=\frac12\bigl(\lambda^2-2\lambda\partial_x\gamma_0-\gamma_0^2\bigr).
  \label{eq:2.53}
\end{equation}
This identity holds on both domains and for both signs of $\lambda$.
\end{proof}

\section{Momentum characterizations of global existence}
\subsection{Momentum characterization on the real line}
\begin{proof}[Proof of Theorem 1.1]
We first assume $\lambda>0$.
\proofpart{Sufficiency.}
Suppose $m_0\ge0$. By \eqref{eq:2.20} and \eqref{eq:2.16},
\begin{equation}
  m(t,q(t,\xi))=\frac{m_0(\xi)}{q_\xi(t,\xi)^2}\ge0,\qquad 0\le t<T^*.
  \label{eq:3.1}
\end{equation}
The flow is onto, so $m(t,\cdot)\ge0$. The reconstruction formulas \eqref{eq:2.21} read
\begin{align}
u(t,x)+\frac\lambda2-u_x(t,x)&=\int_{-\infty}^x e^{y-x}m(t,y)\dd y\ge0,\notag\\
     u(t,x)+\frac\lambda2+u_x(t,x)&=\int_x^\infty e^{x-y}m(t,y)\dd y\ge0, \label{eq:3.2}
\end{align}
and hence, by \eqref{eq:2.7},
\begin{equation}
  |u_x(t,x)|\le u(t,x)+\frac\lambda2\le\frac1{\sqrt2}\norm{u_0}_{H^1(\R)}+\frac\lambda2.
  \label{eq:3.3}
\end{equation}
The continuation criterion \eqref{eq:2.2} therefore implies $T^*=\infty$.
\proofpart{Necessity.}
Suppose $T^*=\infty$. The contrapositive of Lemma 2.5, applied at every initial time, gives
\begin{equation}
  u_x(t,x)\ge-\left|u(t,x)+\frac\lambda2\right|\qquad(t\ge0,\ x\in\R).
  \label{eq:3.4}
\end{equation}
Fix $t\ge0$. The decay of $u(t,\cdot)$ gives
\begin{equation}
  \lim_{x\to-\infty}\left(u(t,x)+\frac\lambda2\right)=\frac\lambda2>0.
  \label{eq:3.5}
\end{equation}
If a zero exists, choose $a$ in the positive tail and let $b>a$ be the first zero to its right. Then
\begin{equation}
  u(t,x)+\frac\lambda2>0\quad(a\le x<b),\qquad u(t,b)+\frac\lambda2=0.
  \label{eq:3.6}
\end{equation}
On $[a,b)$, \eqref{eq:3.4} yields
\begin{align}
\partial_x\!\left[e^x\left(u(t,x)+\frac\lambda2\right)\right]&=e^x\left(u_x(t,x)+u(t,x)+\frac\lambda2\right)\ge0,\notag\\
     u(t,x)+\frac\lambda2&\ge e^{a-x}\left(u(t,a)+\frac\lambda2\right)>0. \label{eq:3.7}
\end{align}
Letting $x\uparrow b$, we obtain the contradiction
\begin{equation}
  0=u(t,b)+\frac\lambda2\ge e^{a-b}\left(u(t,a)+\frac\lambda2\right)>0.
  \label{eq:3.8}
\end{equation}
Therefore
\begin{equation}
  u(t,x)+\frac\lambda2>0\qquad(t\ge0,\ x\in\R).
  \label{eq:3.9}
\end{equation}

Use $F$ from \eqref{eq:2.31}. By \eqref{eq:2.16}, $F(t,\xi)>0$ for $t>0$. For $0\le r\le t$,
\begin{equation}
  q_\xi(r,\xi)>0,\qquad F(t,\xi)-F(r,\xi)=\int_r^t\frac{\dd\tau}{q_\xi(\tau,\xi)}\ge0.
  \label{eq:3.10}
\end{equation}
The momentum identity \eqref{eq:2.33}, positivity \eqref{eq:3.9}, and weights \eqref{eq:3.10} give
\begin{equation}
  m_0(\xi)+2\partial_\xi\frac1{F(t,\xi)}\ge0\qquad(t>0,\ \xi\in\R).
  \label{eq:3.11}
\end{equation}

For any fixed $a<b$, the flow equation and \eqref{eq:2.7} imply
\begin{align}
0<q(r,b)-q(r,a)&=b-a+\int_0^r\!\bigl[u(\tau,q(\tau,b))-u(\tau,q(\tau,a))\bigr]\dd\tau\notag\\
     &\le b-a+\sqrt2\norm{u_0}_{H^1(\R)}r. \label{eq:3.12}
\end{align}
Applying \eqref{eq:2.32}, we obtain
\begin{align}
0<\int_a^b\frac{\dd\xi}{F(t,\xi)}&\le\left(\int_0^t\frac{\dd r}{q(r,b)-q(r,a)}\right)^{-1}\notag\\
     &\le\left(\int_0^t\frac{\dd r}{b-a+\sqrt2\norm{u_0}_{H^1(\R)}r}\right)^{-1}\longrightarrow0\qquad(t\to\infty). \label{eq:3.13}
\end{align}

If $m_0(x_0)<0$, continuity gives a fixed interval $[a,b]$ on which $m_0<0$. For $\xi\in[a,b]$, integration of \eqref{eq:3.11} gives
\begin{align}
  \frac1{F(t,\xi)}&\ge\frac1{F(t,a)}-\frac12\int_a^\xi m_0(y)\dd y\notag\\
    &\ge-\frac12\int_a^\xi m_0(y)\dd y.
  \label{eq:3.14}
\end{align}
Integrating once more,
\begin{align}
\int_a^b\frac{\dd\xi}{F(t,\xi)}&\ge-\frac12\int_a^b\int_a^\xi m_0(y)\dd y\dd\xi\notag\\
     &=-\frac12\int_a^b(b-y)m_0(y)\dd y>0. \label{eq:3.15}
\end{align}
The last quantity is independent of $t$, contradicting \eqref{eq:3.13}. Hence $m_0\ge0$.

For $\lambda<0$, set
\begin{equation}
  v(t,x)=-u(t,-x),\qquad\norm{v(t)}_{H^s}=\norm{u(t)}_{H^s},\qquad s>5/2.
  \label{eq:3.16}
\end{equation}
In distributions,
\begin{align}
&v_t-v_{txx}+3vv_x-\lambda v_x-2v_xv_{xx}-vv_{xxx}\notag\\
     &\qquad=-\bigl(u_t-u_{txx}+3uu_x+\lambda u_x-2u_xu_{xx}-uu_{xxx}\bigr)(t,-x)=0. \label{eq:3.17}
\end{align}
Thus $v$ solves CH with parameter $-\lambda>0$. The reflection is invertible and preserves the Sobolev norms, so the maximal times coincide. Its initial momentum is
\begin{equation}
  v_0(x)-\partial_x^2v_0(x)+\frac{-\lambda}{2}=-m_0(-x).
  \label{eq:3.18}
\end{equation}
The positive-parameter result gives
\begin{equation}
  T^*=\infty\quad\Longleftrightarrow\quad -m_0(-x)\ge0\quad(x\in\R)\quad\Longleftrightarrow\quad\lambda m_0\ge0.
  \label{eq:3.19}
\end{equation}
This proves Theorem 1.1 for every $\lambda\ne0$.
\end{proof}

\subsection{Momentum characterization on the circle}
\begin{proof}[Proof of Theorem 1.2]
Let $\lambda\in\R$.
\proofpart{Sufficiency.}
If $m_0$ has one sign, so does $m(t,\cdot)$ by \eqref{eq:2.20}. Using the periodic lift in \eqref{eq:2.21}, we obtain
\begin{align}
(u+\frac\lambda2)^2-u_x^2&=\left(\int_{-\infty}^x e^{y-x}m(t,y)\dd y\right)\!\left(\int_x^\infty e^{x-y}m(t,y)\dd y\right)\ge0,\notag\\
     |u_x(t,x)|&\le\left|u(t,x)+\frac\lambda2\right|\notag\\ &\le C\norm{u_0}_{H^1(\T)}+\frac{|\lambda|}{2}\qquad(0\le t<T^*), \label{eq:3.20}
\end{align}
The last bound uses \eqref{eq:2.7}. By \eqref{eq:2.2}, $T^*=\infty$.
\proofpart{Necessity.}
Suppose $T^*=\infty$. The contrapositive of Lemma 2.5 gives
\begin{equation}
  u_x(t,x)\ge-\left|u(t,x)+\frac\lambda2\right|\qquad(t\ge0,\ x\in\T).
  \label{eq:3.21}
\end{equation}
Fix $t$ and suppose $u(t,a)+\lambda/2>0$. Before any first zero $b>a$, \eqref{eq:3.21} gives
\begin{align}
\partial_x\!\left[e^x\left(u(t,x)+\frac\lambda2\right)\right]&=e^x\left(u_x(t,x)+u(t,x)+\frac\lambda2\right)\ge0,\notag\\
     u(t,x)+\frac\lambda2&\ge e^{a-x}\left(u(t,a)+\frac\lambda2\right)>0\qquad(a\le x<b). \label{eq:3.22}
\end{align}
Letting $x\uparrow b$ would give
\begin{equation}
  0=u(t,b)+\frac\lambda2\ge e^{a-b}\left(u(t,a)+\frac\lambda2\right)>0.
  \label{eq:3.23}
\end{equation}
There is therefore no zero to the right of $a$. Periodicity gives positivity for every $x$.

Otherwise $u(t,x)+\lambda/2\le0$ everywhere. If equality holds at $a$, \eqref{eq:3.21} gives
\begin{align}
\partial_x\!\left[e^{-x}\left(u(t,x)+\frac\lambda2\right)\right]&=e^{-x}\left(u_x(t,x)-u(t,x)-\frac\lambda2\right)\ge0,\notag\\
     0\ge u(t,x)+\frac\lambda2&\ge e^{x-a}\left(u(t,a)+\frac\lambda2\right)=0\qquad(x\ge a). \label{eq:3.24}
\end{align}
Periodicity then gives $u(t,\cdot)+\lambda/2\equiv0$. At each time, exactly one of the following holds:
\begin{equation}
  u(t,\cdot)+\lambda/2>0,\qquad u(t,\cdot)+\lambda/2<0,\qquad u(t,\cdot)+\lambda/2\equiv0.
  \label{eq:3.25}
\end{equation}
By the mean conservation \eqref{eq:2.8},
\begin{equation}
  \int_\T\left(u(t,x)+\frac\lambda2\right)\dd x=\int_\T\left(u_0(x)+\frac\lambda2\right)\dd x.
  \label{eq:3.26}
\end{equation}
Thus the alternative in \eqref{eq:3.25} is independent of time. In particular,
\begin{equation}
  \int_\T\left(u_0+\frac\lambda2\right)\dd x=0\quad\Longrightarrow\quad u_0\equiv-\frac\lambda2\quad\Longrightarrow\quad m_0\equiv0.
  \label{eq:3.27}
\end{equation}

Suppose the conserved mean is nonzero. Uniqueness of the characteristic equation and periodicity of $u$ give
\begin{equation}
  q(t,\xi+|\T|)=q(t,\xi)+|\T|.
  \label{eq:3.28}
\end{equation}
Hence $q_\xi$ and $F$ are periodic. Applying \eqref{eq:2.32} over one period yields
\begin{align}
0<\int_\T\frac{\dd\xi}{F(t,\xi)}&\le\left(\int_0^t\frac{\dd r}{q(r,a+|\T|)-q(r,a)}\right)^{-1}\notag\\
     &=\left(\int_0^t\frac{\dd r}{|\T|}\right)^{-1}=\frac{|\T|}{t}\longrightarrow0\qquad(t\to\infty). \label{eq:3.29}
\end{align}
In the momentum identity \eqref{eq:2.33}, the weights satisfy $q_\xi>0$ and $F(t,\xi)-F(r,\xi)\ge0$ for $0\le r\le t$, by \eqref{eq:2.16} and \eqref{eq:2.31}.

If the mean in \eqref{eq:3.26} is positive, then $u+\lambda/2>0$ for all times. Thus
\begin{equation}
  m_0(\xi)+2\partial_\xi\frac1{F(t,\xi)}\ge0,\qquad t>0.
  \label{eq:3.30}
\end{equation}
If $m_0$ is negative somewhere, choose fixed $a<b$ within one period with $m_0<0$ on $[a,b]$. Integrating \eqref{eq:3.30} gives
\begin{equation}
  \frac1{F(t,\xi)}\ge\frac1{F(t,a)}-\frac12\int_a^\xi m_0(y)\dd y\ge-\frac12\int_a^\xi m_0(y)\dd y\qquad(a\le\xi\le b).
  \label{eq:3.31}
\end{equation}
Integrating over $[a,b]$, we obtain
\begin{equation}
  \int_\T\frac{\dd\xi}{F(t,\xi)}\ge\int_a^b\frac{\dd\xi}{F(t,\xi)}\ge-\frac12\int_a^b(b-y)m_0(y)\dd y>0,
  \label{eq:3.32}
\end{equation}
contrary to \eqref{eq:3.29}. Therefore $m_0\ge0$.

If the mean is negative, then $u+\lambda/2<0$ for all times, and \eqref{eq:2.33} instead gives
\begin{equation}
  m_0(\xi)+2\partial_\xi\frac1{F(t,\xi)}\le0,\qquad t>0.
  \label{eq:3.33}
\end{equation}
If $m_0$ is positive somewhere, continuity gives a fixed interval $[a,b]$ within one period with $m_0>0$ there. Integration from $\xi$ to $b$ yields
\begin{equation}
  \frac1{F(t,\xi)}\ge\frac1{F(t,b)}+\frac12\int_\xi^b m_0(y)\dd y\ge\frac12\int_\xi^b m_0(y)\dd y.
  \label{eq:3.34}
\end{equation}
Hence
\begin{equation}
  \int_\T\frac{\dd\xi}{F(t,\xi)}\ge\int_a^b\frac{\dd\xi}{F(t,\xi)}\ge\frac12\int_a^b(y-a)m_0(y)\dd y>0,
  \label{eq:3.35}
\end{equation}
again contradicting \eqref{eq:3.29}. Thus $m_0\le0$ in this case. Theorem 1.2 follows.
\end{proof}
\subsection{Momentum characterization for MOCH}
\begin{proof}[Proof of Theorem 1.3]
Let $\lambda\ne0$ and $\gamma_0\in H^s$, $s>3/2$, on $\R$ or $\T$. By Lemma 2.9, the associated CH initial velocity \eqref{eq:2.43} belongs to $H^{s+1}$, $s+1>5/2$, and the two maximal existence times coincide. The momentum is precisely \eqref{eq:1.7}. On $\R$, Theorem 1.1 gives
\begin{equation}
T^*=\infty\quad\Longleftrightarrow\quad \lambda m_0(x)\ge0\qquad(x\in\R).
  \label{eq:3.36}
\end{equation}
On $\T$, Theorem 1.2 and the same lifespan identity give
\begin{equation}
T^*=\infty\quad\Longleftrightarrow\quad m_0\ge0\text{ on }\T\quad\text{or }m_0\le0\text{ on }\T.
  \label{eq:3.37}
\end{equation}
This proves Theorem 1.3 on both domains.
\end{proof}
\begin{proof}[Proof of Remarks 1.1, 1.2 and 1.3]
The energy and amplitude bounds \eqref{eq:1.21} and \eqref{eq:1.25} follow from Lemma 2.3, applied to the CH velocities.

The global derivative bound follows from \eqref{eq:3.3} and \eqref{eq:3.20}. For $\lambda<0$ on $\R$, the reflected solution \eqref{eq:3.16} satisfies
\begin{equation}
  |u_x(t,x)|=|v_x(t,-x)|\le v(t,-x)-\frac\lambda2=\left|u(t,x)+\frac\lambda2\right|.
  \label{eq:3.38}
\end{equation}
The finite-time CH assertion \eqref{eq:1.23} is precisely the continuation criterion \eqref{eq:2.2}.

For MOCH, \eqref{eq:2.44} yields, in the global case,
\begin{align}
\norm{\gamma(t)}_\infty&\le\norm{u_x(t)}_\infty+\norm{\gamma(t)-u_x(t)}_\infty\notag\\
     &\le\norm{u(t)}_\infty+\frac{|\lambda|}{2}+\norm p_2\norm{\gamma_0}_2+\frac{\norm{p_x}_\infty}{2|\lambda|}\norm{\gamma_0}_2^2. \label{eq:3.39}
\end{align}
The right-hand side is uniformly bounded for $t\ge0$. In the finite-time case,
\begin{align}
\left|\inf_x\gamma(t,x)-\inf_xu_x(t,x)\right|&\le\norm{\gamma(t)-u_x(t)}_\infty\notag\\
     &\le\norm p_2\norm{\gamma_0}_2+\frac{\norm{p_x}_\infty}{2|\lambda|}\norm{\gamma_0}_2^2. \label{eq:3.40}
\end{align}
Together with the CH continuation criterion \eqref{eq:2.2}, this gives
\begin{equation}
  \liminf_{t\uparrow T^*}\inf_x\gamma(t,x)=-\infty,
  \label{eq:3.41}
\end{equation}
which proves the finite-time assertion in \eqref{eq:1.28}.
\end{proof}

\section{First breaking and terminal energy}
\begin{proof}[Proof of Theorem 1.4]
Let $u_0\in H^s(\R)$, $s>5/2$, and $T^*<\infty$. We use the flow $q$ from \eqref{eq:2.15} and set
\begin{equation}
  \rho(t,\xi)=\sqrt{q_\xi(t,\xi)},\qquad 0\le t<T^*,\quad \xi\in\R.
  \label{eq:4.1}
\end{equation}
All constants below are independent of $t<T^*$; they may depend on $\lambda$, $u_0$ and $T^*$. We use the spatial regularity in \eqref{eq:2.36}. Here $C_b^j$ denotes bounded continuous functions with bounded continuous derivatives through order $j$.

The conservative theory of Holden and Raynaud \cite[Theorem 2.9 and Proposition 3.16]{r24} gives
\begin{gather}
  u(T^*)\in H^1(\R)\cap C^{1/2}(\R), \label{eq:4.2} \\
  \norm{u(t)-u(T^*)}_\infty\longrightarrow0\qquad(t\uparrow T^*). \label{eq:4.3}
\end{gather}
Young's inequality, \eqref{eq:2.6} and \eqref{eq:2.7} give the bounds used below:
\begin{equation}
  \norm u_\infty+\norm P_\infty+\norm{P_x}_\infty+\norm P_2+\norm{P_x}_2\le C.
  \label{eq:4.4}
\end{equation}
\proofpart{Step 1. Terminal flow and its label regularity.}
We write $(f\circ q)(t,\xi)=f(t,q(t,\xi))$. Differentiating \eqref{eq:1.3} once in space and using $P_{xx}=P-u^2-u_x^2/2-\lambda u$ gives
\begin{equation}
  \frac d{dt}u_x(t,q(t,\xi))=\left(-\frac12u_x^2+u^2+\lambda u-P\right)(t,q(t,\xi)).
  \label{eq:4.5}
\end{equation}
Together with \eqref{eq:2.16}, this gives
\begin{align}
\rho_t&=\frac12(u_x\circ q)\rho,\notag\\
     \rho_{tt}&=\frac12\bigl((u^2+\lambda u-P)\circ q\bigr)\rho, \label{eq:4.6}
\end{align}
with initial values
\begin{equation}
  \rho(0,\xi)=1,\qquad \rho_t(0,\xi)=\frac12 u_{0,x}(\xi).
  \label{eq:4.7}
\end{equation}
The coefficient of this linear ODE is bounded by \eqref{eq:4.4}. Gronwall's inequality gives
\begin{equation}
  \sup_{0\le t<T^*}\bigl(\norm{\rho(t)}_\infty+\norm{\rho_t(t)}_\infty+\norm{\rho_{tt}(t)}_\infty\bigr)<\infty.
  \label{eq:4.8}
\end{equation}

On every $[0,\tau]$, $\tau<T^*$, \eqref{eq:2.36} permits differentiation of the integral form of \eqref{eq:4.6} in $\xi$. In particular,
\begin{align}
 \partial_\xi(u\circ q)&=(u_x\circ q)q_\xi=2\rho\rho_t, \label{eq:4.9}\\
 \partial_\xi\!\left[\frac12(u^2+\lambda u-P)\circ q\right]&=(2u\circ q+\lambda)\rho\rho_t-\frac12(P_x\circ q)\rho^2, \label{eq:4.10}\\
 \rho_{tt\xi}&=\frac12\bigl((u^2+\lambda u-P)\circ q\bigr)\rho_\xi+(2u\circ q+\lambda)\rho^2\rho_t-\frac12(P_x\circ q)\rho^3, \label{eq:4.11}\\
 \rho_\xi(0,\xi)&=0,\qquad \rho_{t\xi}(0,\xi)=\frac12u_{0,xx}(\xi). \label{eq:4.12}
\end{align}
The forcing terms in \eqref{eq:4.11} are bounded by \eqref{eq:4.4} and \eqref{eq:4.8}; $u_{0,xx}$ is bounded by $s>5/2$. A second application of Gronwall's inequality gives
\begin{equation}
  \sup_{0\le t<T^*}\bigl(\norm{\rho_\xi(t)}_\infty+\norm{\rho_{t\xi}(t)}_\infty+\norm{\rho_{tt\xi}(t)}_\infty\bigr)<\infty.
  \label{eq:4.13}
\end{equation}
Thus both $\rho$ and $\rho_t$ are Lipschitz in time with values in $C_b^1(\R)$, and
\begin{align}
\rho(t,\cdot)&\longrightarrow\rho(T^*,\cdot),\notag\\
 \rho_t(t,\cdot)&\longrightarrow\rho_t(T^*,\cdot)\qquad\text{in }C_b^1(\R). \label{eq:4.14}
\end{align}
The subscripts on the terminal quantities denote these limits and the corresponding spatial derivatives.

The bounds on $u$, together with \eqref{eq:4.14}, give
\begin{gather}
  \norm{q(t,\cdot)-q(r,\cdot)}_\infty\le C|t-r|, \label{eq:4.15} \\
  q_\xi(t,\cdot)=\rho(t,\cdot)^2\longrightarrow\rho(T^*,\cdot)^2\quad\text{uniformly}. \label{eq:4.16}
\end{gather}
Define $q(T^*,\xi)=\lim_{t\uparrow T^*}q(t,\xi)$. The convergence is uniform in $\xi$, and
\begin{equation}
  q(T^*,\cdot)\in C^1(\R),\qquad q_\xi(T^*,\xi)=\rho(T^*,\xi)^2\ge0.
  \label{eq:4.17}
\end{equation}
This proves the assertion about $q(T^*)$ in \eqref{eq:1.29}.
\proofpart{Step 2. Breaking labels.}
The displacement bound in \eqref{eq:2.19} passes to the limit; hence $q(T^*,\cdot)$ is onto and $|q(T^*,\xi)-\xi|\le CT^*$. For fixed $t<T^*$, the compact set $\{u_x(r):0\le r\le t\}\subset C_0(\R)$ and the bounded displacement of $q$ imply
\begin{equation}
  \lim_{R\to\infty}\sup_{0\le r\le t}\sup_{|\xi|\ge R}|u_x(r,q(r,\xi))|=0.
  \label{eq:4.18}
\end{equation}
By \eqref{eq:2.16}, $\rho(t,\cdot)-1\in C_0(\R)$. Uniform convergence therefore gives
\begin{equation}
  \rho(T^*,\cdot)-1\in C_0(\R).
  \label{eq:4.19}
\end{equation}
Define the compact set
\begin{equation}
  Z=\{\xi\in\R:\rho(T^*,\xi)=0\}=\{\xi\in\R:q_\xi(T^*,\xi)=0\}.
  \label{eq:4.20}
\end{equation}
If $Z$ were empty, then
\begin{equation}
\inf_{\xi\in\R}\rho(T^*,\xi)>0
  \label{eq:4.21}
\end{equation}
and
\begin{equation}
\limsup_{t\uparrow T^*}\norm{u_x(t)}_\infty=\limsup_{t\uparrow T^*}\left\|\frac{2\rho_t(t)}{\rho(t)}\right\|_\infty<\infty,
  \label{eq:4.22}
\end{equation}
contrary to \eqref{eq:2.2}. Thus $Z\ne\varnothing$.

For $z\in Z$, positivity of $\rho(t,z)$ for $t<T^*$ implies $\rho_t(T^*,z)\le0$. If $\rho_t(T^*,z)=0$, then \eqref{eq:4.6} and its bounded coefficient give
\begin{equation}
  |\rho(t,z)|+|\rho_t(t,z)|\le C\int_t^{T^*}\bigl(|\rho(r,z)|+|\rho_t(r,z)|\bigr)\dd r.
  \label{eq:4.23}
\end{equation}
Backward Gronwall would imply $\rho(0,z)=0$, contrary to $\rho(0,z)=1$. Continuity and compactness therefore give
\begin{equation}
  \rho_t(T^*,z)<0\quad(z\in Z),\qquad \min_{z\in Z}\bigl(-\rho_t(T^*,z)\bigr)>0.
  \label{eq:4.24}
\end{equation}
\proofpart{Step 3. Vanishing momentum at breaking labels.}
For $t<T^*$, differentiating $u_x\circ q=2\rho_t/\rho$ and using \eqref{eq:2.20} give
\begin{align}
  (u_{xx}\circ q)\rho^4&=2\rho\rho_{t\xi}-2\rho_t\rho_\xi, \label{eq:4.25} \\
  m_0&=(u\circ q+\frac\lambda2)\rho^4-2\rho\rho_{t\xi}+2\rho_t\rho_\xi. \label{eq:4.26}
\end{align}
By \eqref{eq:4.3}, \eqref{eq:4.2} and uniform convergence of the flow,
\begin{equation}
  u(t,q(t,\cdot))\longrightarrow u(T^*,q(T^*,\cdot))\quad\text{uniformly}.
  \label{eq:4.27}
\end{equation}
Passing to the limit in \eqref{eq:4.26} therefore yields
\begin{align}
m_0(\xi)&=\left(u(T^*,q(T^*,\xi))+\frac\lambda2\right)\rho(T^*,\xi)^4-2\rho(T^*,\xi)\rho_{t\xi}(T^*,\xi)\notag\\
 &\quad+2\rho_t(T^*,\xi)\rho_\xi(T^*,\xi). \label{eq:4.28}
\end{align}
Since $\rho(T^*,\cdot)$ is nonnegative and $C^1$, every $z\in Z$ satisfies
\begin{equation}
  \rho(T^*,z)=\rho_\xi(T^*,z)=0.
  \label{eq:4.29}
\end{equation}
Equations \eqref{eq:4.28} and \eqref{eq:2.20} now give
\begin{equation}
  m_0(z)=0,\qquad m(t,q(t,z))=\frac{m_0(z)}{q_\xi(t,z)^2}=0.
  \label{eq:4.30}
\end{equation}
This proves \eqref{eq:1.30}.
\proofpart{Step 4. Compression and concentrated energy.}
The conservative energy representation in \cite[Definition 2.5, Theorems 3.10 and 3.13]{r24}, written in the original particle label and using $q_\xi=\rho^2$ and $(u_x\circ q)\rho=2\rho_t$, gives
\begin{equation}
  \norm{u_0}_{H^1}^2-\norm{u(T^*)}_{H^1}^2=4\int_Z\rho_t(T^*,\xi)^2\dd\xi.
  \label{eq:4.31}
\end{equation}
Since $\rho(T^*,z)=0$ for $z\in Z$, \eqref{eq:4.8} and \eqref{eq:4.6} give
\begin{align}
0<\rho(t,z)&=-\int_t^{T^*}\rho_t(r,z)\dd r\le C(T^*-t),\notag\\
     |\rho_{tt}(t,z)|&\le C\rho(t,z)\le C(T^*-t). \label{eq:4.32}
\end{align}
Integrating once and twice, uniformly for $z\in Z$, yields
\begin{align}
|\rho_t(t,z)-\rho_t(T^*,z)|&\le C(T^*-t)^2,\notag\\
     |\rho(t,z)+(T^*-t)\rho_t(T^*,z)|&\le C(T^*-t)^3. \label{eq:4.33}
\end{align}
For $t<T^*$, the flow is an increasing diffeomorphism, so the area formula gives
\begin{equation}
  |q(t,Z)|=\int_Zq_\xi(t,\xi)\dd\xi=\int_Z\rho(t,\xi)^2\dd\xi.
  \label{eq:4.34}
\end{equation}
Combining \eqref{eq:4.33} with the compactness of $Z$, we obtain
\begin{align}
&\left|\frac{4|q(t,Z)|}{(T^*-t)^2}-4\int_Z\rho_t(T^*,\xi)^2\dd\xi\right|\notag\\
     &\qquad\le 4\int_Z\left|\left(\frac{\rho(t,\xi)}{T^*-t}\right)^2-\rho_t(T^*,\xi)^2\right|\dd\xi\le C(T^*-t)^2. \label{eq:4.35}
\end{align}
Equations \eqref{eq:4.31} and \eqref{eq:4.35} prove \eqref{eq:1.31} and complete the proof.
\end{proof}
\begin{proof}[Proof of Remark 1.4]
By Lemma 2.9, $u_0\in H^{s+1}(\R)$, $s>3/2$, and the maximal times coincide. We use $q$, $\rho$ and $Z$ from the preceding proof. The conservative construction of Luo, Qiao and Yin \cite[Theorem 1.3, Corollary 1.4 and Section IV]{r28}, and uniqueness in \cite[Corollary 2.4]{r22}, give
\begin{equation}
  \gamma(t)\rightharpoonup\gamma(T^*)\quad\text{in }L^2(\R)\qquad(t\uparrow T^*).
  \label{eq:4.36}
\end{equation}
\proofpart{Step 1. Weighted terminal limit.}
The MOCH equation \eqref{eq:1.6} and \eqref{eq:4.6} give
\begin{align}
(\gamma\circ q)\rho&=2\rho_t+\bigl((\gamma-u_x)\circ q\bigr)\rho,\notag\\
     \partial_t\bigl[(\gamma\circ q)\rho\bigr]&=\frac12\bigl((\gamma-u_x)\circ q\bigr)(\gamma\circ q)\rho+\lambda(u\circ q)\rho. \label{eq:4.37}
\end{align}
The first identity, \eqref{eq:2.44} and \eqref{eq:4.8} bound the weighted product uniformly. Changing variables by $x=q(t,\xi)$ and using \eqref{eq:2.9}, \eqref{eq:2.6} and \eqref{eq:4.4}, we obtain
\begin{gather}
\norm{(\gamma\circ q)\rho}_{L^2(d\xi)}=\norm{\gamma_0}_2, \notag\\
     \norm{\partial_t[(\gamma\circ q)\rho]}_{L^2(d\xi)}+\norm{\partial_t[(\gamma\circ q)\rho]}_\infty\le C. \label{eq:4.38}
\end{gather}
Thus there is a common limit $\Gamma\in L^2(\R)\cap C_b(\R)$ with
\begin{gather}
(\gamma(t)\circ q(t))\rho(t)\longrightarrow\Gamma\quad\text{in }L^2(d\xi)\text{ and uniformly},\notag\\
     \norm{\Gamma}_2^2=\norm{\gamma_0}_2^2,\qquad \Gamma(\xi)=2\rho_t(T^*,\xi)\quad(\xi\in Z). \label{eq:4.39}
\end{gather}
The last identity uses \eqref{eq:2.44} and $\rho(T^*,\xi)=0$ on $Z$.
\proofpart{Step 2. Equality of the terminal energy defects.}
Monotonicity and \eqref{eq:4.17} give, for $a<b$,
\begin{equation}
  q(T^*,a)=q(T^*,b)\quad\Longrightarrow\quad\int_a^b\rho(T^*,\xi)^2\dd\xi=0\quad\Longrightarrow\quad[a,b]\subset Z.
  \label{eq:4.40}
\end{equation}
Thus $q(T^*)$ maps $\R\setminus Z$ bijectively onto $\R\setminus q(T^*,Z)$, with a $C^1$ inverse on each component. The area formula gives
\begin{equation}
  |q(T^*,Z)|\le\int_Z q_\xi(T^*,\xi)\dd\xi=0.
  \label{eq:4.41}
\end{equation}

On compact subsets of $\R\setminus Z$, \eqref{eq:4.14} and \eqref{eq:4.39} yield local uniform convergence after composition with the inverse flow. The known weak limit \eqref{eq:4.36}, represented continuously on this open set, therefore satisfies
\begin{equation}
  \gamma(T^*,q(T^*,\xi))=\frac{\Gamma(\xi)}{\rho(T^*,\xi)},\qquad\xi\notin Z.
  \label{eq:4.42}
\end{equation}
Changing variables on these components gives
\begin{equation}
  \norm{\gamma(T^*)}_2^2=\int_{\R\setminus Z}|\Gamma(\xi)|^2\dd\xi.
  \label{eq:4.43}
\end{equation}
Combining \eqref{eq:4.39}, \eqref{eq:4.43} and \eqref{eq:4.31}, we obtain
\begin{align}
\norm{\gamma_0}_2^2-\norm{\gamma(T^*)}_2^2&=\int_Z|\Gamma(\xi)|^2\dd\xi\notag\\
     &=4\int_Z\rho_t(T^*,\xi)^2\dd\xi\notag\\
     &=\norm{u_0}_{H^1}^2-\norm{u(T^*)}_{H^1}^2. \label{eq:4.44}
\end{align}
This proves \eqref{eq:1.32}.
\end{proof}

\section{Compression and numerical interpretation}
Throughout this section, $\lambda\ne0$, $u_0\in H^s(\mathbb R)$, $s>5/2$, and
$T^*<\infty$. We use $\rho$ and $Z$ from \eqref{eq:4.1} and \eqref{eq:4.20}. Let $h$ denote
the label mesh spacing, and let $\rho_h$ and $b_h$ approximate $\rho$ and
$2\rho_t$. The numerical integrals use a consistent quadrature rule on
the same label mesh. Write $T_h$ for the first detected zero of $\rho_h$,
and $z_h$ for a label at which that zero is detected. In the diagnostics
below,
\begin{equation}
 \delta=T_h-t. \label{eq:5.1}
\end{equation}
The computed time $T_h$ approximates $T^*$; its discretization error is
distinct from the continuum error in \eqref{eq:4.35}.
\subsection{Compression of the breaking labels}
Theorem 1.4 identifies the rescaled compressed width with the terminal
energy defect in \eqref{eq:1.31}. Its proof also gives the quadratic error bound
\eqref{eq:4.35}.
\clearpage

\noindent
We take $\lambda=1$ and prescribe a smooth initial momentum that vanishes on
$I=[0,1]$. In \eqref{eq:1.6} choose $\gamma_0=\tanh(g/2)$, with $g(x)=x$ on $I$.
Then \eqref{eq:1.7} gives
\[
 m_0(x)=\tfrac12\operatorname{sech}^2(g(x)/2)\bigl(1-g'(x)\bigr).
\]
To specify the extension, put $A(y)=e^{-1/y}$ for $y>0$ and $A(y)=0$
otherwise, and set $S(y)=A(y)/(A(y)+A(1-y))$. For
$r=\operatorname{dist}(x,I)$, take
\[
\begin{aligned}
 g(x)&=S(x+3)\bigl(1-S(x-3)\bigr)\bigl(x+\psi(r)\bigr),\\
 \psi(r)&=\int_0^r A'(y)\bigl(1-S(2y-1)\bigr)\,dy .
\end{aligned}
\]
Thus $m_0=0$ on $I$ and $m_0<0$ on $(1,2)$, so Theorem~1.1 gives
$T^*<\infty$. We evolve the CH velocity from
$u_0=-p_x*\gamma_0-\tfrac12 p*\gamma_0^2$.

On a uniform label grid in $[-20,20]$ we use composite trapezoidal
quadrature for the nonlocal term and integrate the characteristic equations
by DOP853 (relative and absolute tolerances $2\times10^{-9}$ and
$2\times10^{-11}$, maximum time step $0.15$). The event defining $T_h$ is
the first zero of $\min_j\rho_h(t,\xi_j)$ over the entire grid.

The width of the image of $I$ is given by
\begin{equation}
  q(t,1)-q(t,0)=\int_0^1\rho(t,\xi)^2\,d\xi. \label{eq:5.2}
\end{equation}
The numerical width is
\begin{equation}
  W_h(\delta)=\int_0^1\rho_h(T_h-\delta,\xi)^2\,d\xi. \label{eq:5.3}
\end{equation}
With integrals evaluated on the label grid, Figure~\ref{fig:plateau}(a) shows
\begin{equation}
  E_h(\delta)=\left|\frac{4W_h(\delta)}{\delta^2}-m_h\right|,
  \qquad m_h=\int_0^1 b_h(T_h,\xi)^2\,d\xi. \label{eq:5.4}
\end{equation}
If $I=Z$, \eqref{eq:4.31} identifies the integral of $(2\rho_t(T^*,\cdot))^2$ over
$I$ with the terminal energy defect, and \eqref{eq:4.35} gives a quadratic error
bound. We use $I$ as the prescribed zero-momentum interval.

\subsection{Resolution of a compressed interval}
For all labels in $I$ to collapse at first breaking, one needs
\begin{equation}
 \rho(T^*,\xi)=0\qquad (0\leq\xi\leq1). \label{eq:5.5}
\end{equation}
Figure~\ref{fig:plateau}(b) measures the largest residual at the computed event,
\begin{equation}
 R_h=\max_{\xi_j\in I}|\rho_h(T_h,\xi_j)|. \label{eq:5.6}
\end{equation}

\newpage
\begin{center}
\includegraphics[width=\linewidth]{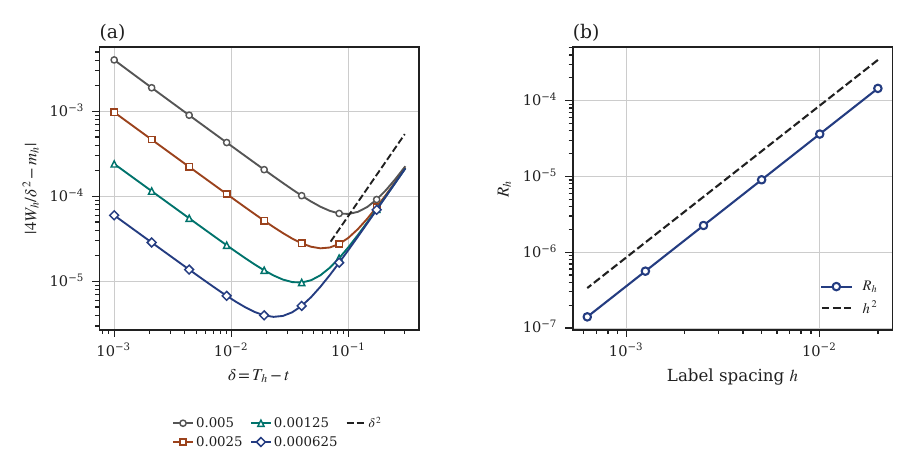}
\end{center}
\refstepcounter{figure}\label{fig:plateau}\noindent Figure \thefigure. Compression of $I=[0,1]$. (a) Error \eqref{eq:5.4} against the time
remaining to the discrete event, for the four finest meshes. (b) Residual
\eqref{eq:5.6} against label spacing, for all six meshes. Dashed lines have slopes
two in $\delta$ and in $h$, respectively.

On the finest grid, a least-squares fit of $\log E_h$ against $\log\delta$
at the six sampled gaps in $[0.1,0.3]$ gives slope $1.985$. At smaller
gaps the curves turn upwards. Even a small terminal residual on $I$ enters
$W_h(\delta)$ before division by $\delta^2$.

The four finest meshes give the following event times, first-zero labels,
and terminal values on $I$:
\begin{center}
\begin{tabular}{rcccc}
\hline
$h$ & $T_h$ & $z_h$ & $m_h$ & $R_h$ \\
\hline
0.005    & 8.77748030 & 1.0550 & 0.04258912 & $9.034\times10^{-6}$ \\
0.0025   & 8.77743760 & 1.0525 & 0.04259152 & $2.253\times10^{-6}$ \\
0.00125  & 8.77742687 & 1.0475 & 0.04259212 & $5.620\times10^{-7}$ \\
0.000625 & 8.77742418 & 1.0450 & 0.04259227 & $1.402\times10^{-7}$ \\
\hline
\end{tabular}
\end{center}
On halving $h$ in the last two rows, $R_h$ decreases by a factor of
$4.008$. This measures spatial resolution at $T_h$. The first discrete zero
lies to the right of $I$ on all four grids. The residuals quantify
compression on $I$ without determining the exact set $Z$.

\clearpage
\noindent\textbf{Conflict of interest.} The authors declare that they have no conflict of interest.

\noindent\textbf{Data availability statement.} The code and data for Figure~\ref{fig:plateau} are included in the accompanying source archive.

\noindent\textbf{Acknowledgments.} Zhaoyang Yin was partially supported by the National Natural Science Foundation of China (Grant No. 12571261).
\phantomsection
\end{document}